\documentstyle[amssymb]{article}

\newcommand{\stopproof}{\hfill \nobreak\medskip $\blacksquare$\\ 
 \hspace*{\fill}} 
\newcommand{\dom}{\mbox{\rm dom}}

\newcommand{\ran}{\mbox{\rm ran}}
\newcommand{\AND}{\mbox{ \rm and }}

\newcommand{\OR}{\mbox{ \rm or }}

\newcommand{\forces}[2]{\Vdash_{#1} \mbox{``} #2 \mbox{''}}
\newcommand{\notforces}[2]{\not\Vdash_{#1} \mbox{``} #2 \mbox{''}}
\newcommand{\proof}{{\bf Proof:} \ }

\newcommand{\PP}{{\Bbb P}}

\newcommand{\DD}{{\Bbb D}}
\newcommand{\RR}{{\Bbb R}}
\newcommand{\CC}{{\Bbb C}}
\newcommand{\KK}{{\Bbb K}}
\newcommand{\QQ}{{\Bbb Q}}

\newcommand{\restricts}{\! \upharpoonright \!}
\newcommand{\lcard}{\, \mid \!}
\newcommand{\rcard}{\! \mid \,}
\newcommand{\betan}{\beta{\Bbb N} \setminus {\Bbb N}}

\newcommand{\pomega}{{\cal P}(\omega)}
\newcommand{\pomegaf}{{\cal P}(\omega)/[\omega]^{<\aleph_0}}
\newcommand{\card}[1]{\lcard #1 \rcard}
\newtheorem{theor}{Theorem}[section]

\newtheorem{defin}{Definition}[section]

\newtheorem{lemma}{Lemma}[section]
\newtheorem{quest}{Question}[section]
\newtheorem{claim}{Claim}

\author{Saharon Shelah and Juris Stepr\={a}ns}
\title{Somewhere Trivial Autohomeomorphisms}
\begin{document}
\thanks{This research was completed while the first author
was supported by the BSF and Rutgers University and the second author
by NSERC and Rutgers University. This is publication Number 427 in the first
author's list of publications.}
    \maketitle
\begin{abstract}
    It is show to be consistent that there is a non-trivial
autohomeomorphism of $\betan$ while all such autohomeomorphisms are
trivial on some open set. The model used is one due to Velickovic in
which, coincidentally, Martin's Axiom also holds.
\end{abstract}
\bibliographystyle{amsplain}
\section{Introduction}
An automorphism of of $\pomegaf$ --- or, equivalently, an
autohomeomorphism of $\betan$ --- is said to be {\em trivial} if there is a
bijection between cofinite subsets of $\omega$ which induces it; an
automorphism is said to be {\em somewhere trivial} if its restriction to
$\cal{P}(A)$ is trivial for some $A\in [\omega]^{\aleph_0}$.
It was shown by Shelah, pages 129 to 152 of \cite{shel.pf}, 
that it is equiconsistent with ZFC that
all automorphisms of $\pomegaf$ are trivial. The
argument which proves this can be viewed as two distinct 
and almost independent arguments. The first part shows that it is consistent
that every
automorphism of $\pomegaf$ is somewhere trivial while the
second part expands on this argument to obtain the consistency of
the assertion that all automorphisms are
indeed trivial. Since the reasoning involved in both parts is, at
least superficially, similar it is natural to ask whether it might
not just be a consequence of the fact that every automorphism is
somewhere trivial, that every automorphism is actually trivial. It is
the purpose of this paper to show that such a theorem does not exist
and hence, the second part of Shelah's argument in \cite{shel.pf} is
indispensable; at the same time this answers Question 205 from \cite{opt.havm}.
 
In order to be more precise the following definitions will be introduced.
\begin{defin}
The relation $\equiv^{\ast}$ has the standard meaning --- namely,
$A\equiv^{\ast} B$ if and only if $\card{A\Delta B}<\aleph_0$ (here,
$A\Delta B = (A\setminus B)\cup (B\setminus A)$). Also.
$A\subseteq^{\ast}B$ is defined to mean that $\card{B\setminus A} < \aleph_0$.
If  $A\subseteq \omega$ then the equivalence class of $A$   with respect to
$\equiv^{\ast}$ will be denoted by $[A]$.
\end{defin}
The notion of triviality can now be precisely formulated.
\begin{defin}
A homomorphism $\Phi :\pomegaf \rightarrow \pomegaf$ is said to be
trivial on $A\subseteq \omega$ if there is $A'\equiv^* A$ and 
a one-to-one function
$f:A'\rightarrow \omega$ such that $\Phi([B]) = [f(B)]$ for every
$B\subseteq A$. A homomorphism will be said to be somewhere trivial if
there is some $A\in [\omega]^{\aleph_0}$ on which it is trivial. A
homomorphism is trivial if it is trivial on $\omega$.
\end{defin} 

It has already been mentioned that it was shown in \cite{shel.pf} that it
is consistent that all automorphisms of $\pomegaf$ are trivial. The
argument relied on the oracle chain condition and it was not clear
what the effect of Martin's Axiom was on the question. This was
partially answered in \cite{step.12} where it was shown that PFA
implies that all  automorphisms of $\pomegaf$ are trivial --- for
related results see \cite{just.orc}. The other half of the answer
was provided by  Velickovic in \nocite{veli.def}\cite{veli.oca} where it is shown that
it is consistent with Martin's Axiom that a nontrivial  automorphism
 of $\pomegaf$ exists.

The following theorem of \cite{veli.oca}  offers an alternate
characterisation of triviality which has proven to be very useful.
\begin{lemma}
    (Velickovic)
If $\Phi :\pomegaf \rightarrow \pomegaf$ is an
automorphism and there exist Borel functions
 $\phi_n$ for $n\in\omega$ and a comeagre set $G\subset \pomega$ such that
 \label{velctble} for every $A\in G$ there is
$n\in\omega$ such that $[\phi_n(A)] = \Phi([A])$ then $\Phi$ is trivial.
\end{lemma}
This is Theorem  2 of \cite{veli.oca} except that in \cite{veli.oca} there is no
reference to the comeagre set $G$; however an inspection of the proof
in \cite{veli.oca}
will reveal that the hypothesis of Theorem 2  can be
weakened to include $G$.
Notice that if $\Phi :\pomegaf \rightarrow \pomegaf$ is a trivial
automorphism then it is simple to find a continuous function
$\phi$ on $\pomega$ such that $[\phi(A)] = \Phi([A])$ for each
$A\subset\omega$. 

The notation  $i_X$ will be used to  denote  the constant 
function whose domain is $X$ and which
has value $i$ at each point in $X$.
Whenever reference is made to a topology on $\pomega$ this
will be to the Cantor set topology under the canonical identification of
$2^{\omega}$ with $\pomega$ --- in other words, a natural 
base for this topology consists of all sets of the form 
$$\{A\subseteq\omega : 1_A\cup 0_{\omega\setminus A}\supseteq g\}$$
where $g$ is a finite partial function from $\omega$ to 2.

The argument to be presented in the next section will be a
modifiction and combination of arguments from pages 129 to 152 of
\cite{shel.pf}, \cite{step.12} and \cite{veli.oca}. For the reader's
benefit, some definitions and lemmas from \cite{shel.pf} will be recalled.
\begin{defin}
    An $\omega_1$-oracle is a sequence $\frak{M} = \{\frak{M}_{\xi} :
\xi\in\omega_1\}$ such that 
\begin{itemize}
    \item $\frak{M}_{\xi}$ is a countable, transitive model of ZFC
without the power set axiom
\item $\xi\in\frak{M}_{\xi}$ and $\frak{M}_{\xi}\models\xi\mbox{ is countable}$
\item $\{\xi\in\omega_1 : A\cap \xi\in\frak{M}_{\xi}\}$ contains a
closed unbounded set
for each $A\subseteq\omega_1$
\end{itemize}
\end{defin}
Notice that the existence of an oracle requires that
$\lozenge_{\omega_1}$ is true.
\begin{defin}
If $\frak{M}$ is an oracle then a partial order $\leq$ on $\omega_1$ 
(or some set coded by $\omega_1$) will be said to satisfy the
$\frak{M}$-chain condition if there is a closed unbounded set $C$
such that  for every $\xi\in C$ and
$A\subseteq \xi$, $A\in \frak{M}_{\xi}$, if $A$ is predense in the  order
$(\xi,\leq\cap(\xi\times\xi))$ then it is predense in 
$(\omega_1,\leq)$.
\end{defin}
Further discussion of these definitions as well as proofs of the following
lemmas can all be found in \cite{shel.pf}.
\begin{lemma} Assume\label{shel2.1} that $\lozenge_{\omega_1}$ holds and $\phi_{\xi}(x)$ is a
$\Pi^1_2$ formula --- possibly with a real parameter --- for each
$\xi\in\omega_1$. Suppose also that there is no $r\in\RR$ such that
$\phi_{\xi}(r)$ holds for all $\xi\in\omega_1$ and that there is still no
such $r$ even after adding a Cohen real. Then there is an oracle
$\frak{M}$ such that any
partial order  $\QQ$  which satisfies the $\frak{M}$-chain condition will not add $r\in\RR$ such that
$\phi_{\xi}(r)$ holds for all $\xi\in\omega_1$.
\end{lemma}
\begin{lemma}
    If  $\{\frak{M}^{\xi} : \xi\in \omega_1\}$ are oracles then there is a
single oracle $\frak{M}$ such that if any partial order satisfies the 
$\frak{M}$-chain condition \label{shel3.1} then it satisfies the $\frak{M}^{\xi}$ chain
condition for 
each $\xi\in \omega_1$. 
\end{lemma}
The oracle $\frak M$ of Lemma~\ref{shel3.1} is easily decribed. It
is the  diagonal union of the oracles
$\{\frak{M}^{\xi} : \xi\in \omega_1\}$. This fact, rather than the
statement of Lemma~\ref{shel3.1},
 will be used in the proof of Lemma~\ref{seqf}.
\begin{lemma}
If $V$  is a model of $\lozenge_{\omega_1}$  then there is, in $V$,
 an oracle $\frak{M}$ \label{2ndcat}
such that if $\QQ$ satisfies the $\frak{M}$-chain condition then
$1\forces{\QQ}{\RR\cap V \mbox{ is second category}}$
\end{lemma}
\begin{lemma} \label{cccorc}
If $\frak{M}$ is  any oracle and $\QQ$ satisfies the 
$\frak{M}$-chain condition then
$\QQ$ satisfies the countable chain condition.
\end{lemma}
\section{The proof}
 The following partial order $\PP$,
was introduced by Velickovic in \cite{veli.oca}
to add a non-trivial automorphism of $\pomegaf$ while doing as 
little else as possible --- at least assuming PFA.
\begin{defin}
    The partial order $\PP$ is defined to consist of all one-to-one functions
$f:A\rightarrow B$ where 
\begin{itemize}
\item $A\subseteq\omega$ and $B\subseteq\omega$ 
\item for all $i\in\omega$ and $n\in\omega$, $f(i)\in (2^{n+1}\setminus 2^n)$ if and only if
 $i\in (2^{n+1}\setminus 2^n)$
  \item $\limsup_{n\rightarrow\omega}\card{(2^{n+1}\setminus 2^n)\setminus
A} = \omega$
and hence, by the previous condition,
$\limsup_{n\rightarrow\omega}\card{(2^{n+1}\setminus 2^n)\setminus
B} = \omega$
\end{itemize} 
The ordering on $\PP$ is $\subseteq^{\ast}$.
\end{defin}

 The terms $2^n$ are not crucial since any sequence of intervals
whose size tends to infinity could equally well have been used.
Further modifications to the partial order are also possible --- some
can be found in \cite{step.18} --- but will not be important in
the present context. It is however, useful to note the following.
\begin{lemma}\label{a} Assume MA$_{\lambda}$. Suppose
that $\eta\leq\lambda$ and that
$$\{f_{\xi} : \xi \in \eta\}$$ is an increasing sequence from
$\PP$. Suppose further that there is $f'$ such that $f' \supseteq^{\ast}
f_{\xi}$ for each $\xi \in \eta$. Then there is $f \in \PP$
such that  $f \supseteq^{\ast}
f_{\xi}$ for each $\xi \in \eta$\end{lemma}
\proof 
It follows from MA$_{\lambda}$ that
there are $A$ and $B$ such that
\begin{itemize}
    \item $A\supseteq^{\ast} \dom({f}_{\xi})$ for each $\xi\in\eta$
    \item $B\supseteq^{\ast} \ran({f}_{\xi})$ for each $\xi\in\eta$
\item  $\limsup_{n\rightarrow\omega}\card{(2^{n+1}\setminus 2^n)\setminus
A} = \omega$
 \item $\limsup_{n\rightarrow\omega}\card{(2^{n+1}\setminus 2^n)\setminus
B} = \omega$
\end{itemize}
 Let 
$f =  f{'}\restricts (A\cap ({f'}^{-1}B))$.
\stopproof
\begin{lemma}
  $\PP$ is countably closed.\label{countably-closed}
\end{lemma}
\proof Given a sequence $\{f_n : n \in \omega\}\subseteq \PP$ 
such that
 $f_n
\subseteq^{\ast} f_{n+1}$ for each $n \in \omega$ choose  inductively  
$k_n$ such that
 $f_{\omega} = \cup\{f_{n}\restricts (\omega \setminus k_{n}) : n 
\in \omega\}$
 is a function.
Now  apply Lemma \ref{a}.
\stopproof
     
>From Lemma \ref{a} it follows that, given a sequence 
$\{f_{\xi} : \xi\in\omega_1\}$, it will be useful
to find an element $f \in \PP$ such that $f_{\xi}\subseteq^{\ast} f$
for each $\xi\in\omega_1$. The following partial order is designed
to do precisely this. 
\begin{defin}
    Given $\{f_{\xi} : \xi\in\mu\}=\frak{F}$ define $\PP(\frak{F})$ 
to be the partial order consisting of all $g\in
\PP$ such that there is some $\xi\in\mu$ such that
$g\equiv^{\ast} f_{\xi}$. The ordering on  $\PP(\frak{F})$ 
is $\subseteq$ as opposed to $\subseteq^{\ast}$ in
$\PP$. 
\end{defin}
\begin{defin}
     For any
$G$ which is  a centred subset of $\PP$ define $\Phi_G :\pomegaf
\rightarrow \pomegaf $ by 
$$\Phi_G([X]) = \left\{ \begin{array}{ll} 
[\{ g(i) : i \in X\}] & \mbox{ if }(\exists g \in G)( X \subseteq \dom(g))\\ 
\left[\omega \setminus \{ g(i) : i \in \omega \setminus X\}\right] &
\mbox{ if }(\exists g \in G)(\omega \setminus X \subseteq \dom(g))
\end{array}\right.$$
If $\Phi$ is a $\PP$-name for an automorphism of $\pomegaf$  then define
$\partial_G{\Phi}(A) = B$ if and only if there is some $p\in  G$ such
that $p\forces{\PP}{\Phi(A) = B}$.
\end{defin}

Velickovic showed that forcing with $\PP$ yields a non-trivial automorphism of $\pomegaf$.
\begin{theor} \label{one-step}
If  $G\subseteq \PP$ is a generic filter on $\PP$ then $\Phi_G$ is
 a non-trivial automorphism of $\pomegaf$.
\end{theor}
\proof 
If it can be shown that $\dom(\Phi_G) = \pomegaf = \ran(\Phi_G)$ then
it is routine to check that $\Phi_G$ 
induces the desired   autohomeomorphism of $\betan$. 
 To see that this is so, assume that $p\in\PP$ and $X \subseteq \omega$
--- since $\PP$ is countably closed, by Lemma \ref{countably-closed},
there is no harm in assuming that
$X\in V$.
It may also be assumed that
$\limsup_{n\rightarrow\omega}\card{(2^{n+1}\setminus 2^n)\setminus
(\dom(p)\cup X)} =\omega$
(otherwise deal with $\omega \setminus X$). It must be shown   that there is
$p' \supseteq p$ such that $p'\forces{\PP}{[X] \in \dom(\Phi_G)}$.
To do
this let $p'\supseteq p$ be any extension satisfying that
 $p'(i)\in (2^{n+1}\setminus 2^n)$ if and only if  $i\in
(2^{n+1}\setminus 2^n)$ for all $i\in X$ and $n\in\omega$. 
A similar proof works for the range of $\Phi$. 
\stopproof
An important fact is the result of
Velickovic \cite{veli.oca}
that if $F$ is $\PP$-generic over a model $V$ of PFA, then
in $V[F]$, not only is there is a non-trivial autohomeomorphism of $\betan$,
but MA also holds.
It will be shown that a closer analysis of this model yields that
in $V[F]$ all autohomeomorphisms of $\betan$ are somewhere trivial.

Loosely speaking, the following theorem will show that 
if $\Phi :\pomegaf \rightarrow \pomegaf$ is a nontrivial
automorphism then it is still nontrivial after adding a Cohen real.
\begin{lemma}
     If $\Phi\in V$ is not trivial and $V'$ is obtained by adding a Cohen
real to \label{cohen} $V$ then, in $V'$, there do not exist Borel functions
$\{\psi_n : n\in\omega\}$ such that for each $C\in\cal{P}(\omega)\cap
V$ there is some $n\in\omega$ such that  $\Phi([C]) = [\psi_n(C)]$.
\end{lemma}
\proof
Suppose that $V'$ is obtained by forcing with the countable partial
order $\CC$ and that $\psi_n$ are $\CC$-names for Borel functions
 such that for each $C\in\cal{P}(\omega)\cap
V$ there is some $n\in\omega$ such that  $\Phi([C]) = [\psi_n(C)]$. Let
$G_n$ be a name for a comeagre set such that $\psi_n\restricts G_n$
is continuous.
 Define $\psi_n^p =\{(A,B) :
p\forces{\CC}{\psi_n(A) = B\AND A\in G_n}\}$. 
Let $D^p_n$ be the closure of the
domain of $\psi_n^p$ and let $E^p_n$ be the closure of the interior
of $D^p_n$ --- note that $D^p_n\setminus E^p_n$ is meagre. 
Let $f_n^p$ be the maximal extension of $\psi_n^p$ to a continuous
function on $E_n^p$. 

It must be that case that the domain of $f_n^p$ is comeagre in $E_n^p$
because if the domain of $f_n^p$ is not comeagre in $E_n^p$ then,
because it is Borel, there must be some open set $U\subset E_n^p$ such
that the set of points in $U$ to which $\psi_n^p$ can be continuously
extended is meagre in $U$.  Since $p\forces{}{\psi_n^p\subset \psi_n}$
and because being a meagre Borel set absolute, it must be that the set
of points in $U$ to which $\psi_n$ can be continuously extended is
also meagre in $U$.  The reason is that the domain $\psi_n^p$ is dense
in $E_n^p$ and so it follows that the domain $\psi_n^p$ is dense in
$U$ and, moreover, not being a point to which a function can be
continuously extended is an absolute property.  This contradicts the
fact that $G_n$ is comeagre.

Now let $M'= \cup\{D^p_n\setminus E^p_n : n\in\omega\AND p\in\CC\}
\cup \{E_n^p\setminus\dom(f_n^p) : n\in \omega\}$ and observe that 
$M'$ is meagre. Now recall the following fact: If $V$ is a model of
{\em ZFC} and $r$ is a Cohen real and $N\in V[r]$ is a meagre set then
there is a meagre set $N'\in V$ such that $N\cap V\subset N'\cap V$.
Let $N$ be a meagre set such that  $ G_n\supset
\pomega\setminus N$ for each $n\in \omega$. Let $M = M'\cup N$.  It is
true in $V'$ that for every $A\in(\pomega\setminus M)\cap V$ there is
some $p\in\CC$ such that $[\psi_n^p(A)] =\Phi([A])$. Since this
statement is arithmetic in the parameters $A$ and $\Phi([A])$ --- and
both of these parameters belong to $V$ --- this must be true in $V$
also. Now apply Lemma
\ref{velctble}. \stopproof

\begin{lemma}
Given $\eta\in\omega_1$, a sequence \label{generic-condition}
$\{f_{\xi} : \xi \in\eta\} = \frak{F_{\eta}}$ and a countable
elementary submodel $\frak{A}\prec (H(\omega_2),\in)$, such that
$\frak{F_{\eta}}\in \frak{A}$, there is $f \in \PP$ which is
$\frak{A}$-generic for $\PP(\frak{F_{\eta}})$. Moreover, for any
extension $\{f_{\xi} : \xi \in \mu\}=\frak{F}_{\mu}$ of
$\frak{F}_{\eta}$ such that $\eta\in\mu\in\omega_1$ and $f_{\eta} =
f$, every $D \in\frak{A}$ is predense in $\PP(\frak{F}_{\mu})$
provided that it is dense in $\PP(\frak{F}_{\eta})$.
\end{lemma}\proof
Let $\{E_k : k\in\omega\}$ enumerate all dense subsets of
$\PP(\frak{F}_{\eta})$ in $\frak{A}$.  Construct sequences $\{g_n :
n\in\omega\}$ and $\{K_n : n\in\omega\}$ such that for all
$n\in\omega$ \begin{itemize} \item $g_n \subseteq g_{n+1}$
\item $K_n < K_{n+1} $
\item $g_{n+1}\restricts K_n = g_n\restricts K_n$
\item there is some $i$ such that $2^{i+1} = K_n$ and
$\card{(2^{i+1}\setminus 2^i)\setminus \dom(g_n)}\geq n$
\item for each bijection $t:K_n\rightarrow K_n$ there $h\in \cap_{j\in n+1}
E_j$ such that $t\cup g_n\restricts(\omega\setminus K_n) \supseteq h$
\end{itemize}   It is easy to see that
this can be done. Hence, it is possible to define $f=\cup\{g_n :
n\in\omega\}$. Notice
 that $\{g : g\supseteq^* f_{\xi}\}$ is dense in
$\PP(\frak{F}_{\eta})$ and definable in $\frak{A}$ --- hence
$f\supseteq^*f_{\xi}$ for each $\xi\in\eta$.
 To check that $f$ has the desired properties suppose
that $g\in\PP(\frak{F}_{\mu})$ for some $\mu\geq \eta$ and that
$f_{\eta} = f$. If $E$ is dense in $\PP(\frak{F}_{\eta})$ then there
is some $m$ such that $E\in \{E_j : j\in m+1\}$ and 
$g\restricts(\omega\setminus K_m)\supseteq
f\restricts (\omega\setminus K_m)$. By extending $g$ if necessary,
it may, without loss of generality, be assumed that
 $g\restricts K_m = t$  and  that $t:K_m\rightarrow K_m$ is a bijection.
It follows that $t\cup g_m\restricts(\omega\setminus K_m) \supseteq h$
for some $h\in E$ and hence $g\supseteq h\in E$.
\stopproof
\begin{lemma}  Suppose that $V$ is a model of $2^{\aleph_0} = \aleph_1$. 
    If  $\Phi$ is \label{seq2} a $\PP$-name for a nowhere trivial automorphism of
$\pomegaf$ and $f\in\PP$ then there is a sequence
$\frak{F}=\{f_{\xi} : \xi \in\omega_1\}\subset \PP$ such that $f_0 = f$ and
$\partial_{\frak{F}}\Phi$ is nowhere trivial.
\end{lemma}
\proof
Let $\{A_{\xi} :
\xi\in\omega_1\}$ be an enumeration of $\pomegaf$ in $V$ and
$$\{\Psi_{\xi} : \cal{P}(C_{\xi}) 
 \rightarrow\cal{P}(B_{\xi}) : \xi\in\omega_1\}$$ enumerate all
possible names for continuous functions
from a Borel comeagre subset of some $\cal{P}(C)$ to some
$\cal{P}(B)$ so that each name occurs cofinally often. 
It suffices to construct $\frak{F}=\{f_{\xi} :
\xi\in\omega_1\}\subseteq \PP$ by induction so that
for every limit ordinal  $\xi$ the following conditions are satisfied
\begin{itemize}
    \item $f_{\xi + n}$ decides, in $\PP$, the values of 
$\Phi(A_{\xi + n -2})$ and
 $\Phi^{-1}(A_{\xi + n -2})$ for $n\geq 2$
\item there is some $C\subset A_{\xi}$ such that 
$1\forces{\PP(\frak{F})}{\partial_{\frak{F}}\Phi(C) \neq [\Psi_{\xi}(C)]}$
\end{itemize}
It is possible to construct $\frak{F}$ inductively because
a failure would mean that for some $\xi\in\omega_1$  it must be the case that
$$f_{\xi}\forces{\PP}{\Phi\mbox{ is trivial on } A_{\xi}}$$ contradicting
that $\Phi$ is a name for a nowhere trivial automorphism of $\pomegaf$.
\stopproof
\begin{lemma} Suppose that $V$ is a model of $\lozenge$ and that
 $\Phi$ is a $\PP$-name for a nowhere trivial automorphism of $\pomegaf$.
Then there is  a sequence\label{seqf}
$\frak{F}=\{f_{\xi} : \xi \in\omega_1\}\subset \PP$ such that
\begin{itemize}
\item $\PP(\frak{F})$ satisfies the countable 
chain condition
\item $\RR\cap V$ is of second category after forcing with
$\PP(\frak{F})$
\item for every $G\subset \PP(\frak{F})$ which
is generic over $V$, for every 
$A\in V\cap\pomega$, $B\in V\cap\pomega$ and  
$\PP(\frak{F})$-names $\Psi_n$ such that for
each $n\in\omega$ 
$$1\forces{\PP(\frak{F})}{ \Psi_n:\cal{P}(A)\rightarrow\cal{P}(B)\mbox{ is continuous}}$$
there is some $C\in  V$  such that
 $\partial_{\frak{F}}\Phi([C])\neq\Psi_n([C])$ for all $n\in\omega$ 
\end{itemize}    
In the last clause the possibility that $C\notin\dom(\Psi_n)$ is
allowed in the sense that if  $C\notin\dom(\Psi_n)$ then 
$\partial_{\frak{F}}\Phi([C])\neq\Psi_n([C])$.
\end{lemma}
\proof
The proof will be rely on constructing a particular oracle which will
guarantee that the three clauses are all satisfied. The only wrinkle
is that the oracle and the sequence $\{f_{\xi} : \xi \in\omega_1\}$
must be constructed simultaneously. The sequence 
$\{f_{\xi} : \xi \in\omega_1\}$ will be obtained by diagonalizing
across $\aleph_1$ such sequences.
 
In particular, let $\frak N$ be any oracle such that forcing with an
$\frak N$-oracle chain condition partial order preserves the fact that
$\RR\cap V$ is of second category --- such an oracle exists by 
Lemma~\ref{2ndcat}. Then construct sequences
$\frak{F}^{\mu}=\{f_{\xi}^{\mu} : \xi \in\omega_1\}\subset \PP$
and $\frak{M}^{\mu}=\{\frak{M}_{\xi}^{\mu} : \xi \in\omega_1\}\subset \PP$  for
$\mu\in\omega_1$  such that\begin{itemize}
			       \item[\bf{a.}] $f^{\xi}_{\mu} = f_{\mu}^{\mu}$
if $\mu < \xi$
\item[\bf{b.}] $\partial_{\frak{F}^{\mu}}\Phi$ is nowhere trivial for $\mu\in\omega_1$
\item[\bf{c.}]for $\mu\in\omega_1$, if $\QQ$ satisfies the $\frak{M}^{\mu }$-chain condition
and $G$ is $\QQ$-generic over $V$ then, in $V[G]$,  for every 
$A\in V\cap\pomega$, $B\in V\cap\pomega$ there do not exist
 $\{\Psi_n : n\in\omega\}$ such that  
 $\Psi_n:\cal{P}(A)\rightarrow\cal{P}(B)$ is continuous and
for all $C\in \cal{P}(A)\cap V$  there exists $n\in\omega$ such that
 $\partial_{\frak{F}^{\mu}}\Phi([C]) = \Psi_n([C])$
\item[\bf{d.}] $\{f_{\mu}^{\xi} : \{\mu,\xi\}\in [\eta]^2\}\in
\frak{M}_{\eta}^{\eta}$ for each $\eta\in\omega_1$
\item[\bf{e.}] $f_{\mu}^{\mu + 1}$ is $\PP(\{f_{\xi}^{\xi} :
\xi\in\mu\})$-generic over $\frak{M}_{\mu + 1}^{\mu}$
\item[\bf{f.}]  $\frak{M}_{\xi}^{\mu}  \in \frak{M}_{\xi}^{\eta}$ if
$\xi\in\mu\in\eta$
\item[\bf{g.}] $\frak{M}^0 = \frak N$
			   \end{itemize}

To see that this suffices let $\frak{M}_{\xi}  = \frak{M}_{\xi + 1}^{\xi}$
and let  $\frak{F} = \{f_{\xi}  = f_{\xi + 1}^{\xi + 1} : \xi\in\omega_1\}$. 
It follows from the remark following  Lemma~\ref{shel3.1} 
that $\{\frak{M}_{\xi} :{\xi}\in\omega_1\}$ is an
oracle and that any partial order which satisfies the
 $\{\frak{M}_{\xi} :{\xi}\in\omega_1\}$-chain condition also satisfies
each
of the $\frak{M}^{\mu}$-chain conditions for $\mu\in\omega_1$. 
Since $f_{\mu}^{\mu + 1}$ is $\PP(\{f_{\xi}^{\xi} :
\xi\in\mu\})$-generic over $\frak{M}_{\mu + 1}^{\mu}$
it follows that
 $\PP(\{f_{\xi} :
\xi\in\omega_1\})$ satisfies the $\{\frak{M}_{\xi} : \xi\in\omega_1\}$
chain condition.
In
particular, this partial order satisfies the $\frak{M}^0$-chain
condition and hence the second clause of the theorem will be
satisfied. That the first clause is satisfied follows from Lemma~\ref{cccorc}.
 So it only remains to be shown that the last clause
is satisfied.

To this end, suppose that $A\in V\cap\pomega$, $B\in V\cap\pomega$ and  
$\PP(\frak{F})$-names $\Psi_n$ are given such that
for each $n\in\omega$ 
$$1\forces{\PP(\frak{F})}{ \Psi_n:
\cal{P}(A)\rightarrow\cal{P}(B)\mbox{ is continuous}}$$
Since $\PP(\frak{F})$ 
satisfies the countable chain condition,
there is some $\gamma\in\omega_1$ such that $\frak{M}_{\gamma}$
models that
for each $n\in\omega$ 
$$1\forces{\PP(\{f_{\xi} : \xi \in\gamma\})}{ \Psi_n:
\cal{P}(A)\rightarrow\cal{P}(B)\mbox{ is continuous}}$$
It now follows that this statement about $\frak{M}_{\gamma} = 
\frak{M}_{\gamma}^{\gamma}$ must be true in $\frak{M}_{\gamma +
1}^{\gamma}$ because $\frak{M}^{\gamma}_{\gamma} 
\in\frak{M}^{\gamma}_{\gamma + 1}$. But it now follows from the fact that
$f_{\gamma}^{\gamma + 1}$ is generic over $\frak{M}^{\gamma}_{\gamma +
1}$ that $$
f_{\gamma}^{\gamma + 1}\forces{\PP(\{f_{\xi}^{\xi} : \xi\in\gamma\})}
{(\exists C\in \cal{P}(A)\cap V)(\forall n\in\omega)
 \partial_{\frak{F}^{\gamma}}\Phi([C]) \neq \Psi_n([C])}$$
Since $f_{\gamma}^{\gamma + 1} \supset^* f^{\gamma+1}_{\mu} = f_{\mu}^{\mu}$ 
for each $\mu\in\omega_1\cap\frak{M}_{\gamma+1}^{\gamma}$  
it follows that $\partial_{\frak{F}}\Phi\restricts\frak{M}_{\gamma +
1}^{\gamma} = \partial_{{\frak{F}}^{\gamma}}\Phi\restricts\frak{M}_{\gamma + 1}$
and hence
$$
f_{\gamma}^{\gamma + 1}\forces{\PP(\{f_{\xi}^{\xi} : \xi\in\gamma\})}
{(\exists C\in \cal{P}(A)\cap V)(\forall n\in\omega)
 \partial_{\frak{F}}\Phi([C]) \neq \Psi_n([C])}$$
Since the necessary dense sets are definable in $\frak{M}_{\gamma +
1}^{\gamma} = \frak{M}_{\gamma}$ it follows that
$$
f_{\gamma}\forces{\PP(\frak{F})}
{(\exists C\in \cal{P}(A)\cap V)(\forall n\in\omega)
 \partial_{\frak{F}}\Phi([C]) \neq \Psi_n([C])}$$ which is what is required.

All that remains to be done is to show that the inductive
construction can be completed. For this, suppose that
$\frak{F}^{\mu}=\{f_{\xi}^{\mu} : \xi \in\omega_1\}$
and $\frak{M}^{\mu}=\{\frak{M}_{\xi}^{\mu} : \xi
\in\omega_1\}$  have been constructed for $\mu\in\eta$.  

If $\eta$ is a limit then it is easy to use Lemma \ref{seq2} in
order to satisfy conditions (a) and (b). If $\eta $ is a successor
then Lemma \ref{generic-condition}
 must also be used in order to satisfy condition (e). To construct $\frak{M}^{\eta}_{\xi}$ for 
$\xi \in\omega_1$ use Lemmas \ref{shel2.1} and Lemma \ref{cohen} to
satisfy condition (c). It is then easy to enlarge the terms of the oracle
to satisfy conditions (d) and (f).
 \stopproof
The proof of the main theorem will require the following definition, which
is a reformulated form of the partial order which appeared in \cite{shel.pf} on
page 134.
\begin{defin}
Given a sequence $\{(W_{\xi},V_{\xi}) : \xi\in\eta\}$ 
define $\QQ(\{(W_{\xi},V_{\xi}) : \xi\in\eta\})$ to be the
partial order which consists of all functions $g$ such that there is
$\Gamma\in [\eta]^{<\aleph_0}$ such that
$$g\equiv^{\ast} \cup\{0_{V_{\xi}}\cup 1_{W_{\xi}\setminus V_{\xi}}
: \xi\in\Gamma\}$$
 The ordering 
on  $\QQ(\{(W_{\xi},V_{\xi}) : \xi\in\eta\})$ is inclusion.
\end{defin}
\begin{theor}
It is consistent, relative to the consistency of ZFC and PFA, 
that\label{mainth} all automorphisms
of $\pomegaf$ are somewhere trivial
but there is, nevertheless, a non-trivial automorphism of $\pomegaf$.
\end{theor}
\proof
It follows from Theorem \ref{one-step} that 
forcing with $\PP$ yields a non-trivial
automorphism of $\pomegaf$.
Hence all that needs to be shown is that
in the resulting model all automorphisms of $\pomegaf$ are somewhere trivial.
To do this suppose that $V$ is a model of PFA and that in this model
 $$1\forces{\PP}{\Phi
\mbox{ is a nowhere trivial automorphism of }\pomegaf}.$$

 Let $G$ be a $V$-generic filter on a countably closed 
partial order --- the Levy collapse of $2^{\aleph_0}$ to $\aleph_1$
for example --- which forces the existence of a $\lozenge$-sequence.
Let
$\frak{F} = \{f_{\xi} : \xi\in\omega_1\}$ be some fixed sequence with the
properties guaranteed by Lemma \ref{seqf}
In particular, Lemma~\ref{seqf} guarantees that
 for every 
$A\in V[G]\cap\pomega$, $B\in V[G]\cap\pomega$ and for every
collection of  
$\PP(\frak{F})$-names $\{\Psi_n : n\in\omega\}$ 
such that, for each $n\in\omega$,
$$1\forces{\PP(\frak{F})}{ \Psi_n:\cal{P}(A)\rightarrow\cal{P}(B)
\mbox{ is continuous}}$$
there is some $C\in \cal{P}(A)\cap V[G] = \cal{P}(A)\cap V$ such that
 $\partial_{\frak{F}}\Phi([C])\neq\Psi_n([C])$ for every $n\in\omega$.
Let $H$ be  $V[G]$-generic for the partial order  $\PP(\frak{F})$. 

Let $\frak M$ be an arbitrary oracle in $V$.
A sequence $\{(W_{\xi},V_{\xi}) : \xi\in\omega_1\}$ will be
constructed in $V[G][H]$ so that
\begin{itemize}
\item  if $\QQ_{\xi} = \QQ(\{(W_{\eta},V_{\eta}) : \eta\in\xi\})$
then $\QQ_{\omega_1}$ satisfies the $\frak{M}$-chain condition
    \item $(W_{\xi},V_{\xi})\in V[G] = V\cap (\pomega)^2$
\item $V_{\xi}\subseteq W_{\xi}\subseteq \omega$
\item $\card{W_{\xi}\cap W_{\eta}} <\aleph_0$ if $\eta\neq\xi$
\item  for  each
$p\in\QQ_{\xi}$ and $\QQ_{\xi}$-name, $Y\in\frak{M}_{\xi}$, for a
subset of $\omega$
$$p\cup 1_{V_{\xi}}\cup 0_{W_{\xi}\setminus
V_{\xi}}\forces{\QQ_{\xi + 1}}{\partial_{\frak F}\Phi([W_{\xi}])\cap [Y] \neq
\partial_{\frak F}\Phi([V_{\xi}])} $$
\item the dense subsets of $\QQ_{\xi + 1}$ which guarantee that the
previous statement is true are predense in $\QQ_{\omega_1}$ 
\end{itemize}

Before continuing, define $\Phi^*(A)\subseteq  \omega$ arbitrarily to satisfy
that $[\Phi^*(A)] = \partial_{\frak F}\Phi([A])$ for each
$[A]\in\dom(\partial_{\frak F}\Phi)$. Next, choose  an almost disjoint family
 $\{W'_{\xi} : \xi\in\omega_1\}$ in the model $V[G]$.
The set $W_{\xi}$ will be chosen so that, among other things,
$W_{\xi}\subseteq W'_{\xi}$ --- this will, of course, guarantee that the
resulting family is almost disjoint.
If this construction succeeds then it is possible to 
proceed as in \cite{step.12} to prove that forcing with 
$\PP(\frak{F})\ast
\QQ_{\omega_1}$ adds a set to which
the partial automorphism $\partial_{\frak F} \Phi$ can not be extended. 

In particular, if $H_1\ast H_2$ is
$\PP(\frak{F})\ast
\QQ_{\omega_1}$ generic then,
  setting $X = \cup\{f^{-1}(\{1\}) : f\in H_2\}$, it 
follows  that $X\cap W_{\xi}\equiv^{\ast}
V_{\xi}$ for each $\xi\in\omega_1$ but, in $V[G][H_1\ast H_2]$, for every
$Y\subseteq\omega$ there is $\beta\in\omega_1$ such that
$\Phi^*(W_{\xi})\cap Y\not\equiv^{\ast}\Phi^*(V_{\xi})$
 for each $\xi\geq\beta$. Just as in \cite{step.12}, it is possible to 
 define a relation
$R$ on $\omega_1$ by $R(\xi,\eta)$ holds if and only if
either $(\Phi^*(W_{\xi})\setminus \Phi^*(V_{\xi}))
\cap \Phi^*(V_{\eta})\neq\emptyset$ or
 $(\Phi^*(W_{\eta})\setminus 
\Phi^*(V_{\eta})\cap \Phi^*(V_{\xi}))\neq\emptyset$. 
It is easy to see that this is a semiopen relation --- as defined in 
\cite{ab.ru.sh} --- 
and that moreover, there is no 
$S\in [\omega_1]^{\aleph_1}$ such that $[S]^2\cap R =\emptyset$. The reason for the
last statement is that otherwise, letting $Y = \cup\{\Phi^*(V_{\xi})
 : \xi \in S\}$ 
would yield a contradiction to the fact that
$\Phi^*(W_{\xi})\cap Y\not\equiv^{\ast}\Phi^*(V_{\xi})$ for all but
countably  many $\xi$.
Hence, by the results of \cite{ab.ru.sh},
there is a proper partial order $\KK$ which adds a 
set $S\in [\omega_1]^{\aleph_1}$ such that $[S]^2\subseteq R$.
This makes the fact that $\partial_{\frak F}\Phi$ 
can not be extended to the set $X$ absolute.
The reason for {\em this} is that if there is a set $Y$ such that $\partial_{\frak F}\Phi([X]) $
can be defined to be $[Y]$, then it must be the case that
$Y\cap\Phi^*(W_{\xi})\equiv^* \Phi^*(V_{\xi})$ for each $\xi\in S$.
But then there is an uncountable set $S'\subseteq S$, as well as $J\in \omega$,
such that  $Y\cap\Phi^*(W_{\xi})\setminus J =
\Phi^*(V_{\xi})\setminus J$ 
for each $\xi\in S'$. It follows that $\Phi^*(V_{\xi})\setminus
J\subseteq Y$ and that  $(\Phi^*(W_{\xi})\setminus\Phi^*(V_{\xi}))\setminus
J\subseteq \omega\setminus Y$ for each $\xi\in S'$. Choosing $\xi$ and
$\zeta$ in $S'$ such that  $\Phi^*(V_{\xi})\cap  J
= \Phi^*(V_{\zeta})\cap  J$ and  $\Phi^*(W_{\xi})\cap  J
= \Phi^*(W_{\zeta})\cap  J$ yields the desired contradiciton.

The iteration $\DD\ast\PP(\frak{F})\ast\QQ_{\omega_1}\ast\KK$ is proper and only
$\aleph_1$ dense sets in it need be met in order to obtain 
$S$ and the set $X$ such that $\partial_{\frak F}\Phi$ can not
be extended to include $[X]$ in its domain. Let $f_{\omega_1}$ be the element of $\PP$ obtained by
forcing with $\PP(\{frak{F})$ and Lemma \ref{a}
 and note that, in $V$,
$f_{\omega_1}\forces{\PP}{\Phi\mbox{ does not extend to }X}$ because
$f_{\omega_1}\forces{\PP}{\partial_{\frak F}\Phi\subseteq \Phi}$

Hence it may be assumed that the construction breaks down at some point
$\mu\in\omega_1$. What can go wrong?
First, there are certain predense sets required at stage $\mu$ which 
must remain predense in the partial order $\QQ_{\mu + 1}$. It is shown
on page 134 of \cite{shel.pf} 
that, for each predense set $E\subseteq \QQ_{\mu }$,  
there is a dense
open set $\cal{W}\subseteq\cal{P}(W'_{\mu})$ such that if $W\in\cal{W}$ 
then, letting $W_{\mu} = W$, $E$ remains predense in
$\QQ_{\mu + 1}$ for {\em any} $V=V_{\mu}\subseteq W_{\mu}$. Note
that  $\cal{W}$ is closed under  the operation of taking infinite
subsets. Recall that $\frak F$ was chosen so that
 $\cal{P}(W'_{\mu})\cap V[G]$ is 
of second category in $\cal{P}(W'_{\mu})$ in the model $V[G][H_1]$.
It follows that 
it may be assumed that $W_{\mu}\in V[G]$ and that $W_{\mu} \in \cal {O}$ for
every open set $\cal {O}\subseteq\cal{P}(W_{\mu}')$ 
which is definable from $\frak{M}_{\mu}$ and 
$\{(W_{\xi},V_{\xi}) : \xi\in\mu \}$ --- but note that
$\QQ_{\mu}$ is definable from $\{(W_{\xi},V_{\xi}) : \xi\in\mu \}$.
Hence the only possible problem is
 that it is not possible to find $V_{\mu}\subseteq W_{\mu}$
satisfying the required properties --- namely, there is no
$V_{\mu}\subseteq W_{\mu}$ such that
for  each
$p\in\QQ_{\mu}$ and $\QQ_{\mu}$-name, $Y\in\frak{M}_{\mu}$, for a
subset of $\omega$
$$p\cup 1_{V_{\mu}}\cup 0_{W_{\mu}\setminus
V_{\mu}}\forces{\QQ_{\mu + 1}}{\partial_{\frak F}\Phi([W_{\mu}])\cap [Y] 
\not\equiv^{\ast}
\partial_{\frak F}\Phi([V_{\mu}])} .$$

To see that this can not happen it will
be necessary to discuss forcing in $\QQ_{\mu + 1}$ before $\QQ_{\mu +
1}$ has been defined --- namely, before $V_{\mu}$ has been defined.
This will be done by defining $\Vdash_{\ast}{}$ as follows: If
$p\in\QQ_{\mu}$, $X$ and $Y$ are $\QQ_{\mu}$ names for subsets of $\omega$
 and $f:W_{\mu}\rightarrow 2$ is a partial function, 
define $$(p, f)\forces{\ast}{X \not\equiv^{\ast} Y}$$
if and only if for each $p'\in\QQ_{\mu}$ such that $p'\cup p\cup f$
is a function and for each $n\in\omega$ there is $p{''}$ such that
\begin{itemize}
    \item  $p{''}\cup p'\cup p\cup f$ is a function
\item $p{''}\forces{\QQ_{\mu}}{k\in X\Delta Y}$ for some $k\geq n$
\end{itemize}
\begin{claim} If the following conditions are satisfied
\begin{itemize}
\item $(p, f)\forces{\ast}{X'\cap  X\not\equiv^{\ast}  Y}$
\item 
$X'$ is a $\QQ_{\mu}$ name belonging 
to $\frak{M}_{\mu}$ 
\item $X$ and $Y$  are subsets of $\omega$ in $V[G]$ 
\item $V_{\mu} = f^{-1}(\{1\})$
\end{itemize} then\label{one}
 then $p\cup f\forces{\QQ_{\mu +
1}}{X\cap X'\not\equiv^{\ast}Y}$.\end{claim}
The way to see this is to note that the following statements are all equivalent
\begin{itemize}
\item $p{''}\forces{\QQ_{\mu}}{k\in X\cap X'\Delta Y}$
\item  $k\in X\setminus Y$ and $p{''}\forces{\QQ_{\mu}}{k\in  X'}$ or
$k\in Y\setminus X$ and $p{''}\forces{\QQ_{\mu}}{k\notin   X'}$
\item  $k\in X\setminus Y$ and $p{''}\forces{\QQ_{\mu + 1}}{k\in  X'}$ or
$k\in Y\setminus X$ and $p{''}\forces{\QQ_{\mu + 1}}{k\notin   X'}$ so long as the name
$X'$ is still a $\QQ_{\mu + 1}$ name; in other words, the
antichains in $\QQ_{\mu}$ 
deciding membership in $X'$  remain maximal in $\QQ_{\mu + 1}$.
\end{itemize}
This last equivalence is guaranteed by the choice of $W_{\mu}$,
 because the relevant dense 
sets are definable from $X'$, which belongs to $\frak{M}_{\mu}$,
 and $\QQ_{\mu}$. 
The claim now follows from the definition  of $\Vdash_{\ast}{}$.
\medskip 

 It will be now be shown that it is possible to choose
$V_{\mu}\subseteq W_{\mu}$  such that 
for  each
$p\in\QQ_{\mu}$ and $\QQ_{\mu}$-name, $Y\in\frak{M}_{\mu}$, for a
subset of $\omega$
$$p\cup 1_{V_{\mu}}\cup 0_{W_{\mu}\setminus
V_{\mu}}\forces{\QQ_{\mu + 1}}{\Phi^*(W_{\mu})\cap Y 
\not\equiv^{\ast}\Phi^*(V_{\mu})} $$
If this is not possible then it follows from Claim 1 that there 
is no $V_{\mu}$ such that
$V_{\mu}\subseteq W_{\mu}$ and such that for  each
$p\in\QQ_{\mu}$ and each $\QQ_{\mu}$-name for a subset of $\omega$,
 $Y\in\frak{M}_{\mu}$
$$(p, 1_{V_{\mu}}\cup 0_{W_{\mu}\setminus
V_{\mu}})\forces{\ast}{\partial_{\frak F}\Phi([W_{\mu}])
\cap [Y] \not\equiv^{\ast}
\partial_{\frak F}\Phi([V_{\mu}])} . $$
It will be shown that this implies
that there is are
continuous functions $\Psi_n$, for $n\in\omega$, such that
for each $C\subset W_{\mu}$ there is $n\in\omega$ such that
 $\Psi_n(C)\equiv^{\ast}\Phi^*(C)$, thus 
 contradicting
the fact that $\frak{M}$ is being assumed to satisfy the
 conclusion of Lemma \ref{seqf}. 

Here is how to conclude this. 
For each 
$p\in\QQ_{\mu}$ and $A\subset W_{\mu}$ define $q_{p}(A) = (p,
 1_A \cup 0_{W_{\mu}\setminus A})$ and
$\bar{q}_p(A) = p\cup  1_A \cup 0_{W_{\mu}\setminus A}$.
Given $p$ and $r$ in $\QQ_{\mu}$, $n\in\omega$ 
 and $Y\in\frak{M}_{\mu}$, a $\QQ_{\mu}$-name
 for a subset of $\omega$,
 define $\psi_{p,r,n,Y}(A)$ to be the set
$$ \{k\in\Phi^*(W_{\mu})\setminus n : (\forall p')(p'\cup p\cup \bar{q_r}(A)
\mbox{ is not a function }\OR
p'\forces{\QQ_{\mu}}{k\in Y}\}$$
 for $A\in\dom(\psi_{p,r,n,Y}) = \{A\subset W_{\mu} : \bar{q}_r(A)\cup p\mbox{ 
is a function }\}$.

It will first be shown that 
for each $A\subset W_{\mu}$ there are $p$, $r$, $n$ and $Y$ such that
$\Phi^*(A) \equiv^{\ast}
\psi_{p,r,n,Y}(A)$. To see see this recall that
there is some $r_A\in \QQ_{\mu}$ and some $Y_A\in\frak{M}_{\mu}$, a $\QQ_{\mu}$-name
 for a subset of $\omega$, 
such that
 $$q_{r_A}(A)
\forces{\ast}{\Phi^*(W_{\mu})\cap
Y_A\not\equiv^{\ast}\Phi^*(A)}$$
fails to be true and $q_{r_A}(A)$ is a function. Hence,
there is some $p_A\in \QQ_{\mu}$
such that $p_A\cup \bar{q}_{r_A}(A)$ is a function and there is some
$n_A\in\omega$ such that for each $p'$ and  $k\geq n_A$ either
$p'\cup p_A\cup \bar{q}_{r_A}(A)$ is not a function or
$p'\notforces{\QQ_{\mu}}{k\in
(\Phi^*(W_{\mu})\cap
Y_{A})\Delta(\Phi^*(A))}$. 

 If $k \in \psi_{p_A,r_A,n_A,Y_A}(A)\setminus \Phi^*(A)$
and $k > n_A$ then the definition of $\psi_{p_A,r_A,n_A,Y_A}(A)$ implies 
that $p'\forces{\QQ_{\mu}}{{k\in Y_{A}}}$
whenever $p'$ is such that $p'\cup p_A\cup \bar{q}_{r_A}(A)$ 
is a function. But
this means that 
$p'\forces{\QQ_{\mu}}{k\in
(\Phi^*(W_{\mu})\cap
Y_{A})\Delta(\Phi^*(A))}$ contradicting the choices of 
 $p_A$, $r_A$, $ n_A$ and $Y_A$.
 
On the other hand, suppose that
 $k$ belongs to $(\Phi^*(A)\cap\Phi^*(W_{\mu}))\setminus n_A$. 
If $k\not\in \psi_{p_A,r_A,n_A,Y_A}(A)$
then there exists $p'$ such that 
$p'\cup p_A\cup \bar{q}_{r_A}(A)$ is a function and
$p'\notforces{\QQ_{\mu}}{k\in Y_A}$. It follows that there is 
$p'' \supset p'$ such that $p''\forces{\QQ_{\mu}}
{k\not\in Y_A}$. If $p''\cup p_A \cup \bar{q}_{r_A}(A)$ is a function then
$p''\forces{\QQ_{\mu}}{k\in
(\Phi^*(W_{\mu})\cap
Y_{A})\Delta(\Phi^*(A))}$ once again
 contradicting the choice of $ p_A$, $r_A$, $ n_A$ and $ Y_A$.
But why should  $p''\cup p_A \cup \bar{q}_{r_A}(A)$ be a function?

The fact that it is possible to choose 
$p{''}$ such that $p{''}\cup p_A \cup \bar{q}_{r_A}(A)$ 
is a function follows
from the choice of $W_{\mu}$. Recall that $W_{\mu}$ was chosen so
that $W_{\mu}\in\cal{O}$ for every open set 
$\cal{O}\subseteq \cal{P}(W_{\mu}')$ which is definable from
$\frak{M}_{\mu}$ and $\QQ_{\mu}$.
Moreover, the set $\cal{O}$ consisting of all $W\subseteq
W_{\mu}'$ such that there exists $e\in [W]^{<\aleph_0}$ such that for all
$\epsilon :e\rightarrow 2$ there exists
$p_{\epsilon}$ satisfying one of the following three conditions
\begin{itemize} 
\item $\epsilon\cup p_A\cup r_A$ is not a function 
\item $ p_{\epsilon}\forces{\QQ_{\mu}}{k\notin Y_{A}}$ and
$p_{\epsilon}\cup p_A \cup r_A\cup \epsilon$ is a function and
$p_{\epsilon}\restricts  W = \epsilon$
\item $\epsilon\cup p_A \cup r_A \forces{\QQ_{\mu}}{k\in Y_A}$
\end{itemize}
 is easily seen to be dense and open in
$W_{\mu}'$ and to be definable from $\QQ_{\mu}$ and $Y_A$, $p_A$ and
$r_A$ --- all of which belong to $\frak{M}_{\mu}$.  Hence $W_{\mu}$ 
belongs to this dense open set $\cal{O}$. Now let $e\in
[W_{\mu}]^{<\aleph_0}$ witness this fact; in other words, letting
$\epsilon = (p_A\cup r_A \cup1_A\cup 0_{W_\mu\setminus
A})\restricts  e$,
 there is
some $p_{\epsilon}$ such that
\begin{itemize}
    \item $p_{\epsilon}\forces{\QQ_{\mu}}{k\notin Y_{A}}$
\item $p_{\epsilon}\cup p_A \cup r_A$ is a function 
\item $p_{\epsilon}\restricts  W_{\mu} = p_A
\cup r_A \cup 1_A\cup 0_{W_\mu\setminus
A})\restricts  e$
\end{itemize}
because the other two alternatives are not possible in light of the
fact that $k\not\in\psi_{p_A,r_A,n_A,Y_A}(A)$ and
 $A\in\dom(\psi_{p_A,r_A,n_A,Y_A})$.
It follows that setting $p{''}= p_{\epsilon}$ yields the desired
condition.

Now observe that the functions $\psi_{p,r,n,Y}$ are
all Borel. This
 contradiction to Lemma \ref{seqf} finishes the proof of the theorem.

\stopproof
\section{Remarks and Open Questions}
It is worth noting that not only has it been shown that it is
consistent that there is a nontrivial automorphism yet all
automorphisms are somewhere trivial, but also that this is consistent
with MA$_{\omega_1}$. Combining the arguments of this paper with
those of \cite{step.18} it is possible to show that that it is
consistent (even with MA$_{\omega_1}$) that every autohomeomorphism
of $\betan$ is somewhere trivial while any two P-points have the same
topological type in the sense that there is an autohomeomorphism of
$\betan$ mapping one to the other. In this model it will of course
follow that there are $2^{\frak{c}}$ autohomeomorphisms of $\betan$
which raises the following question.
\begin{quest}
    Is it consistent that there are only $2^{\aleph_0}$ automorphisms
of $\pomegaf$ but that  there is, nevertheless, a nontrivial automorphism?
\end{quest}

The nontrivial automorphism constructed with Velickovic's order is
much more than somewhere trivial --- the collection of subsets of
$\omega$ where it is trivial forms a maximal ideal. Given any
automorphism $\Phi :\pomegaf\rightarrow\pomegaf$ define
$\cal{J}(\Phi)$ to be the collection of sets on which $\Phi$ is
trivial. It is not difficult to check that $\cal{J}(\Phi)$ is always
an ideal but it is not clear what else can be said about it. 
\begin{quest}
    Does MA$_{\omega_1}$ imply that $\cal{J}(\Phi)\neq
[\omega]^{<\aleph_0}$ for each autmorphism $\Phi$ of $\pomegaf$? 
\end{quest}
\begin{quest}
    Does MA$_{\omega_1}$ imply 
 that $\cal{J}(\Phi)$ is 
the intersection of maximal ideals
 for every
automorphism $\Phi$ of $\pomegaf$? 
\end{quest}
\makeatletter \renewcommand{\@biblabel}[1]{\hfill#1.}\makeatother


\begin{thebibliography}{1}

\bibitem{ab.ru.sh}
U.~Abraham, M.~Rubin, and S.~Shelah, {\em On the consistency of some partition
  theorems for continuous colorings and the structure of $\aleph_1$-dense real
  order types}, Annals of Pure and Applied Logic {\bf 29} (1985), 123--206.

\bibitem{opt.havm}
K.~P. Hart and J.~van Mill, {\em Open problem on $\beta\omega$}, Open Problems
  in Topology (G.~M. Reed and J.~van Mill, eds.), North-Holland, pp.~99--125.

\bibitem{just.orc}
W.~Just, {\em A modification of {S}helah's oracle-c.c. with applications},
  Transactions of the American Mathematical Society {\bf ???} (1989), ??--??

\bibitem{shel.pf}
S.~Shelah, {\em Proper forcing}, Lecture Notes in Mathematics, vol. 940,
  Springer-Verlag, Berlin, 1982.

\bibitem{step.12}
S.~Shelah and J.~Stepr\={a}ns, {\em {P}{F}{A} implies all automorphism are
  trivial}, Proc. Amer. Math. Soc. {\bf 104} (1988), 1220--1225.

\bibitem{step.18}
J.~Stepr\={a}ns, {\em Martin's {A}xiom and the transitivity of {P}-points},
  Submitted to Israel Journal of Mathematics, 19.

\bibitem{veli.def}
B.~Velickovic, {\em Definable automorphisms of $\cal{P}(\omega)/fin$}, Proc.
  Amer. Math. Soc {\bf 96} (1986), 130--135.

\bibitem{veli.oca}
B.~Velickovic, {\em {O}{C}{A} and automorphisms of $\cal{P}(\omega)/fin$}, preprint,
  1990.

\end{thebibliography}
\end{document}